# Exponential Ratio-Product Type Estimators Under Second Order Approximation In Stratified Random Sampling

[1]Rajesh Singh, [1]Prayas Sharma and [2]Florentin Smarandache

[1]Department of Statistics, Banaras Hindu University

Varanasi-221005, India

[2]Chair of Department of Mathematics, University of New Mexico, Gallup, USA

**Abstract**

Singh et al. (20009) introduced a family of exponential ratio and product type estimators in stratified random sampling. Under stratified random sampling without replacement scheme, the expressions of bias and mean square error (MSE) of Singh et al. (2009) and some other estimators, up to the first- and second-order approximations are derived. Also, the theoretical findings are supported by a numerical example.

## 1. INTRODUCTION

In survey sampling, it is well established that the use of auxiliary information results in substantial gain in efficiency over the estimators which do not use such information. However, in planning surveys, the stratified sampling has often proved needful in improving the precision of estimates over simple random sampling. Assume that the population U consist of L strata as $U=U_1, U_2,\ldots,U_L$ . Here the size of the stratum $U_h$ is $N_h$, and the size of simple random sample in stratum $U_h$ is $n_h$, where h=1, 2,---,L.

When the population mean of the auxiliary variable, $\overline{X}$, is known, Singh et al. (2009) suggested a combined exponential ratio-type estimator for estimating the population mean of the study variable $(\overline{Y})$ :

$$t_{1S} = \bar{y}\exp\left[\frac{\bar{X} - \bar{x}_{st}}{\bar{X} + \bar{x}_{st}}\right] \tag{1.1}$$

where,

$$\bar{y}_h = \frac{1}{n_h}\sum_{i=1}^{n_h} y_{hi}, \qquad \bar{x}_h = \frac{1}{n_h}\sum_{i=1}^{n_h} x_{hi},$$

$$\bar{y}_{st} = \sum_{h=1}^{L} w_h \bar{y}_h, \quad \bar{x}_{st} = \sum_{h=1}^{L} w_h \bar{x}_h, \text{ and } \quad \bar{X} = \sum_{h=1}^{L} w_h \bar{X}_h.$$

The exponential product-type estimator under stratified random sampling is given by

$$t_{2S} = \bar{y}\exp\left[\frac{\bar{X} - \bar{x}_{st}}{\bar{x}_{st} + \bar{X}}\right] \tag{1.2}$$

Following Srivastava (1967) an estimator $t_{3s}$ in stratified random sampling is defined as :

$$t_{3S} = \bar{y}\exp\left[\frac{\bar{X} - \bar{x}_{st}}{\bar{x}_{st} + \bar{X}}\right]^{\alpha} \tag{1.3}$$

where $\alpha$ is a constant suitably chosen by minimizing MSE of $t_{3S}$. For $\alpha=1$, $t_{3S}$ is same as conventional exponential ratio-type estimator whereas for $\alpha = -1$, it becomes conventional exponential product type estimator.

Singh et al. (2008) introduced an estimator which is linear combination of exponential ratio-type and exponential product-type estimator for estimating the population mean of the study variable $(\bar{Y})$ in simple random sampling. Adapting Singh et al. (2008) estimator in stratified random sampling we propose an estimator $t_{4s}$ as :

$$t_{4S} = \bar{y}\left[\theta\exp\left[\frac{\bar{X} - \bar{x}_{st}}{\bar{x}_{st} + \bar{X}}\right] + (1-\theta)\exp\left[\frac{\bar{X} - \bar{x}_{st}}{\bar{x}_{st} + \bar{X}}\right]\right] \tag{1.4}$$

where $\theta$ is the constant and suitably chosen by minimizing mean square error of the estimator $t_{4S}$. It is observed that the estimators considered here are equally efficient when terms up to first order of approximation are taken. Hossain et al. (2006) and Singh and Smarandache (2013) studied some estimators in SRSWOR under second order approximation. Koyuncu and Kadilar (2009, 2010) ), have studied some estimators in stratified random sampling under second order approximation. To have more clear picture about the best estimator, in this study we have derived the expressions of MSE's of the estimators considered in this paper up to second order of approximation in stratified random sampling.

## 3. Notations used

Let us define, $e_0 = \frac{\bar{y}_{st} - \bar{y}}{\bar{y}}$ and $e_1 = \frac{\bar{x}_{st} - \bar{x}}{\bar{x}}$,

such that

$$E(e_0) = E(e_1) = E(e_2) = 0,$$

$$V_{rs} = \sum_{h=1}^{L} W_h^{r+s} E\left[(\bar{x}_h - \bar{X}_h)^r (\bar{y}_h - \bar{Y}_h)^s\right]$$

To obtain the bias and MSE of the proposed estimators, we use the following notations in the rest of the article:

$$\bar{y}_{st} = \sum_{h=1}^{L} w_h \bar{y}_h = \bar{Y}(1 + e_0),$$

$$\bar{x}_{st} = \sum_{h=1}^{L} w_h \bar{x}_h = \bar{X}(1 + e_1),$$

where $\bar{y}_h$ and $\bar{Y}_h$ are the sample and population means of the study variable in the stratum h, respectively. Similar expressions for X and Z can also be defined.

Also, we have

$$E(e_0^2) = \frac{\sum_{h=1}^{L} w_h^2 \gamma_h S_{yh}^2}{\bar{Y}^2} = V_{200} \,,$$

$$E(e_1^2) = \frac{\sum_{h=1}^{L} w_h^2 \gamma_h S_{xh}^2}{\bar{X}^2} = V_{020},$$

$$E(e_0 e_1) = \frac{\sum_{h=1}^{L} w_h^2 \gamma_h S_{xyh}}{\bar{X}\bar{Y}} = V_{110},$$

where

$$S_{yh}^2 = \frac{\sum_{i=1}^{N_h} (\bar{y}_h - \bar{Y}_h)^2}{N_h - 1}, \qquad S_{xh}^2 = \frac{\sum_{i=1}^{N_h} (\bar{x}_h - \bar{X}_h)^2}{N_h - 1}, \quad S_{zh}^2 = \frac{\sum_{i=1}^{N_h} (\bar{z}_h - \bar{Z}_h)^2}{N_h - 1},$$

$$S_{xyh} = \frac{\sum_{i=1}^{N_h} (\bar{x}_h - \bar{X}_h)(\bar{y}_h - \bar{Y}_h)}{N_h - 1}, \qquad S_{yzh} = \frac{\sum_{i=1}^{N_h} (\bar{y}_h - \bar{Y}_h)(\bar{z}_h - \bar{Z}_h)}{N_h - 1},$$

$$\gamma_h = \frac{1 - f_h}{n_h}, \qquad f_h = \frac{n_h}{N_h}, \qquad w_h = \frac{N_h}{n_h}.$$

Some additional notations for second order approximation:

$$V_{rs} = \sum_{h=1}^{L} W_h^{r+s} \frac{1}{\bar{Y}^r \bar{X}^s} E\left[ \left(\bar{y}_h - \bar{Y}_h\right)^s \left(\bar{x}_h - \bar{X}_h\right)^r \right]$$

where , $C_{rs(h)} = \frac{1}{N_h} \sum_{i=1}^{N_h} \left[ \left(\bar{y}_h - \bar{Y}_h\right)^s \left(\bar{x}_h - \bar{X}_h\right)^r \right],$

$$V_{12} = \sum_{h=1}^{L} W_h^3 \frac{k_{1(h)} C_{12(h)}}{\bar{Y}\bar{X}^2}, \qquad V_{21} = \sum_{h=1}^{L} W_h^3 \frac{k_{1(h)} C_{21(h)}}{\bar{Y}^2 \bar{X}}, \qquad V_{30} = \sum_{h=1}^{L} W_h^3 \frac{k_{1(h)} C_{30(h)}}{\bar{Y}^3},$$

$$V_{03} = \sum_{h=1}^{L} W_h^3 \frac{k_{1(h)} C_{03(h)}}{\bar{X}^3}, \qquad V_{13} = \sum_{h=1}^{L} W_h^4 \frac{k_{2(h)} C_{13(h)} + 3k_{3(h)} C_{01(h)} C_{02(h)}}{\bar{Y}\bar{X}^3},$$

$$V_{04} = \sum_{h=1}^{L} W_h^4 \frac{k_{2(h)} C_{04(h)} + 3k_{3(h)} C_{02(h)}^2}{\bar{X}^4}, \qquad V_{22} = \sum_{h=1}^{L} W_h^4 \frac{k_{2(h)} C_{22(h)} + k_{3(h)} \left(C_{01(h)} C_{02(h)} + 2C_{11(h)}^2\right)}{\bar{Y}^2 \bar{X}^2},$$

where $k_{1(h)} = \dfrac{(N_h - n_h)(N_h - 2n_h)}{n^2 (N_h - 1)(N_h - 2)}$,

$$k_{2(h)} = \frac{(N_h - n_h)(N_h + 1)N_h - 6n_h (N_h - n_h)}{n^3 (N_h - 1)(N_h - 2)(N_h - 3)},$$

$$k_{3(h)} = \frac{(N_h - n_h)N_h (N_h - n_h - 1)(n_h - 1)}{n^3 (N_h - 1)(N_h - 2)(N_h - 3)}.$$

## 4. First Order Biases and Mean Squared Errors under stratified random sampling

The expressions for biases and MSE,s of the estimators $t_{1S}$, $t_{2S}$ and $t_{3S}$ respectively, are :

$$\text{Bias}(t_{1S}) = \bar{Y}\left[\frac{3}{8}V_{02} - \frac{1}{2}V_{11}\right] \tag{4.1}$$

$$\text{MSE}(t_{1S}) = \bar{Y}^2\left[V_{20} + \frac{1}{2}V_{02} - V_{11}\right] \tag{4.2}$$

$$\text{Bias}(t_{2S}) = \bar{Y}\left[\frac{1}{2}V_{11} - \frac{1}{8}V_{02}\right] \tag{4.3}$$

$$\text{MSE}(t_{2S}) = \bar{Y}^2\left[V_{20} + \frac{1}{4}V_{02} + V_{11}\right] \tag{4.4}$$

$$\text{Bias}(t_{3S}) = \overline{Y}\left[\alpha\frac{1}{4}V_{02} + \alpha^2\frac{1}{8}V_{02} - \frac{1}{2}\alpha V_{11}\right] \tag{4.5}$$

$$\text{MSE}(t_{3S}) = \overline{Y}^2\left[V_{20} + \frac{1}{4}\alpha^2 V_{02} - \alpha V_{11}\right] \tag{4.6}$$

By minimizing $\text{MSE}(t_{3S})$, the optimum value of $\alpha$ is obtained as $\alpha_o = \frac{2V_{11}}{V_{02}}$. By putting this optimum value of $\alpha$ in equation (4.5) and (4.6), we get the minimum value for bias and MSE of the estimator $t_{3S}$.

The expression for the bias and MSE of $t_{4s}$ to the first order of approximation are given respectively, as

$$\text{Bias}(t_{4s}) = \overline{Y}\left[\theta\left\{\left\{\frac{3}{8}V_{02} - \frac{1}{2}V_{11}\right\} + (1-\theta)\left\{\frac{1}{2}V_{11} - \frac{1}{8}V_{02}\right\}\right]\right. \tag{4.7}$$

$$\text{MSE}(t_{4S}) = \overline{Y}^2\left[V_{20} + \left(\frac{1}{2} - \theta\right)^2 V_{02} + 2\left(\frac{1}{2} - \theta\right)V_{11}\right] \tag{4.8}$$

By minimizing $\text{MSE}(t_{4S})$, the optimum value of $\theta$ is obtained as $\theta_o = \frac{V_{11}}{V_{02}} + \frac{1}{2}$. By putting this optimum value of $\alpha$ in equation (4.7) and (4.8) we get the minimum value for bias and MSE of the estimator $t_{3S}$. We observe that for the optimum cases the biases of the estimators $t_{3S}$ and $t_{4S}$ are different but the MSE of $t_{3S}$ and $t_{4S}$ are same. It is also observed that the MSE's of the estimators $t_{3S}$ and $t_{4S}$ are always less than the MSE's of the estimators $t_{1S}$ and $t_{2S}$. This prompted us to study the estimators $t_{3S}$ and $t_{4S}$ under second order approximation.

## 5. Second Order Biases and Mean Squared Errors in stratified random sampling

Expressing estimator $t_i$'s(i=1,2,3,4) in terms of $e_i$'s (i=0,1), we get

$$t_{1s} = \bar{Y}(1+e_0)\exp\left[\frac{-e_1}{2+e_1}\right]$$

Or

$$t_{1s} - \bar{Y} = \bar{Y}\left\{e_0 - \frac{e_1}{2} - \frac{1}{2}e_0e_1 + \frac{3}{8}e_1^{\ 2} + \frac{3}{8}e_0e_1^{\ 2} - \frac{7}{48}e^3 - \frac{7}{48}e_0e_1^{\ 3} + \frac{25}{384}e^4\right\} \tag{5.1}$$

Taking expectations, we get the bias of the estimator $t_{1s}$ up to the second order of approximation as

$$Bias_2(t_{1s}) == \frac{\bar{Y}}{2}\left[\left[-V_{11} + \frac{3}{4}V_{02} + \frac{3}{4}V_{12} - \frac{7}{24}V_{03} - \frac{7}{24}V_{13} + \frac{25}{192}V_{04}\right]\right. \tag{5.2}$$

Squaring equation (5.1) and taking expectations and using lemmas we get MSE of $t_{1s}$ up to second order of approximation as

$$MSE(t_{1S}) = E\left[\bar{Y}\left(e_0 - \frac{e_1}{2} + \frac{3}{8}e_1^{\ 2} - \frac{1}{2}e_0e_1 + \frac{3}{8}e_0e_1^{\ 2} - \frac{7}{48}e^3\right)\right]^2$$

Or,

$$MSE(t_{1s}) = \bar{Y}^2E\left[\left\{e_0^{\ 2} + \frac{1}{4}e_1^{\ 2} - e_0e_1 + e_0^{\ 2}e_1^{\ 2} - e_0^{\ 2}e_1 - \frac{3}{8}e_1^{\ 3} - \frac{25}{24}e_0e_1^{\ 3} + \frac{5}{4}e_0e_1^{\ 2} + \frac{55}{192}e_1^{\ 4}\right\}\right] \tag{5.3}$$

Or,

$$MSE_2(t_{1s}) = \bar{Y}^2\left[V_{20} + \frac{1}{4}V_{02} - V_{11} + V_{22} - V_{21} + \frac{5}{4}V_{12} - \frac{25}{24}V_{13} + \frac{55}{192}V_{04}\right] \tag{5.4}$$

Similarly we get the biases and MSE's of the estimators $t_{2S}$, $t_{3S}$ and $t_{4S}$ up to second order of approximation respectively, as

$$\text{Bias}_2(t_{2s}) = \frac{\overline{Y}}{2}\left[V_{11} - \frac{1}{4}V_{02} - \frac{1}{4}V_{12} - \frac{5}{24}V_{13} + \frac{1}{192}V_{04} - \frac{5}{24}V_{03}\right] \tag{5.5}$$

$$\text{MSE}_2(t_{2S}) = \overline{Y}^2\left[V_{20} + \frac{1}{4}V_{02} + V_{11} + \frac{23}{192}V_{04} - \frac{1}{8}V_{03} + \frac{1}{4}V_{12} - \frac{1}{24}V_{13} + V_{21}\right) \tag{5.6}$$

$$\text{Bias}_2(t_{3S}) = \overline{Y}\left[\left(\frac{\alpha^2}{8} + \frac{\alpha}{4}\right)V_{02} + \left(\frac{\alpha^2}{8} + \frac{\alpha}{4}\right)V_{12} - \frac{\alpha}{2}V_{11}\left(\frac{\alpha^2}{8} + \frac{\alpha^3}{48}\right)V_{03} - \left(\frac{\alpha^2}{8} + \frac{\alpha^3}{48}\right)V_{13}\right.$$

$$\left. + \left(\frac{\alpha^2}{32} + \frac{\alpha^3}{32} + \frac{\alpha^4}{384}\right)V_{04}\right] \tag{5.7}$$

$$\text{MSE}_2(t_{3S}) = \overline{Y}^2\left[V_{20} + \frac{\alpha^2}{4}V_{02} - \alpha V_{11} + \left(\frac{\alpha}{2} + \frac{\alpha^2}{2}\right)V_{22} - \alpha V_{21} + \left(\frac{\alpha}{2} + \frac{\alpha^2}{2}\right)V_{22} + \left(\frac{\alpha}{2} + \frac{3\alpha^2}{4}\right)V_{12}\right.$$

$$\left. - \left(\frac{\alpha^2}{4} + \frac{\alpha^2}{8}\right)V_{03} - \left(\frac{3\alpha^2}{4} + \frac{7\alpha^3}{24}\right)V_{13} + \left(\frac{\alpha^2}{16} + \frac{\alpha^3}{16} + \frac{7\alpha^4}{192}\right)V_{04}\right] \tag{5.8}$$

$$\text{Bias}_2(t_{4S}) = E(t_{4S} - \overline{Y}) = \overline{Y}\left[\left(\frac{1}{2} - \alpha\right)V_{11} - \frac{1}{2}\left(\frac{1}{4} - \alpha\right)\{V_{02} + V_{12}\} + \left\{\frac{1}{16}\left(\frac{1}{24} + \alpha\right)\right\}V_{04}\right.$$

$$\left. - \frac{1}{48}(2\alpha + 5)\{V_{03} + V_{13}\}\right] \tag{5.9}$$

$$\text{MSE}_2(t_{4S}) = \overline{Y}^2\left[V_{20} + \left(\frac{1}{2} - \theta\right)^2 V_{02} + \left\{\left(\frac{1}{2} - \theta\right)^2 + \frac{(4\theta - 1)}{4}\right\}V_{22} + \left\{\left(\frac{1}{2} - \theta\right)^2 + \frac{(4\theta - 1)}{4}\right\}V_{12}\right.$$

$$+ \left\{\frac{1}{64}(4\theta - 1)^2 - \frac{1}{24}\left(\frac{1}{2} - \theta\right)(2\theta + 5)\right\}V_{04} + 2\left(\frac{1}{2} - \theta\right)V_{21} + \frac{1}{4}\left(\frac{1}{2} - \theta\right)(4\theta - 1)V_{03}$$

$$-\left\{-\frac{1}{24}(2\theta+5)+\frac{1}{2}\left(\frac{1}{2}-\theta\right)(4\theta-1)\right\}V_{13}\Bigg] \tag{5.10}$$

The optimum value of $\alpha$ we get by minimizing $MSE_2(t_{3S})$. But theoretically the determination of the optimum value of $\alpha$ is very difficult, we have calculated the optimum value by using numerical techniques. Similarly the optimum value of $\theta$ which minimizes the MSE of the estimator $t_{4s}$ is obtained by using numerical techniques.

## 6. Numerical Illustration

For the one natural population data, we shall calculate the bias and the mean square error of the estimator and compare Bias and MSE for the first and second order of approximation.

**Data Set-1**

To illustrate the performance of above estimators , we have considered the natural data given in Singh and Chaudhary (1986, p.162).

The data were collected in a pilot survey for estimating the extent of cultivation and production of fresh fruits in three districts of Uttar- Pradesh in the year 1976-1977.

**Table 6.1**: **Bias and MSE of estimators**

| Estimator | Bias | | MSE | |
|---|---|---|---|---|
| | First order | Second order | First order | Second order |
| $t_{1s}$ | -1.532898612 | -1.475625158 | 2305.736643 | 2308.748272 |
| $t_{2s}$ | 8.496498176 | 8.407682289 | 23556.67462 | 23676.94086 |
| $t_{3s}$ | -1.532898612 | -1.763431841 | 704.04528 | 705.377712 |
| $t_{4s}$ | -5.14408 | -5.0089 | 704.04528 | 707.798567 |

## 7. CONCLUSION

In the Table 6.1 the bias and MSE of the estimators $t_{1S}$, $t_{2S}$, $t_{3S}$ and $t_{4S}$ are written under first order and second order of approximation. The estimator $t_{2S}$ is exponential product-type estimator and it is considered in case of negative correlation. So the bias and mean squared error for this estimator is more than the other estimators considered here. For the classical exponential ratio-type estimator, it is observed that the biases and the mean squared errors increased for second order. The estimator $t_{3S}$ and $t_{4S}$ have the same mean squared error for the first order but the mean squared error of $t_{3S}$ is less than $t_{4S}$ for the second order. So, on the basis of the given data set we conclude that the estimator $t_{3S}$ is best followed by the estimator $t_{4S}$ among the estimators considered here.